\documentclass[12pt]{amsart}
\usepackage{amsmath,amssymb,color,epsfig,graphicx}
\newtheorem{thm}{Theorem}[section]
\newtheorem{cor}[thm]{Corollary}
\newtheorem{lem}[thm]{Lemma}

\theoremstyle{remark}
\newtheorem{rem}{Remark}[section]
\newtheorem{definition}{Definition}[section]

\theoremstyle{definition}

\numberwithin{equation}{section}
\numberwithin{figure}{section}

\newenvironment{pf}{\begin{proof}}{\end{proof}}

\font\nt=cmr7

\def\note#1
{\marginpar
{\nt $\leftarrow$
\par
\hfuzz=20pt \hbadness=9000 \hyphenpenalty=-100 \exhyphenpenalty=-100
\pretolerance=-1 \tolerance=9999 \doublehyphendemerits=-100000
\finalhyphendemerits=-100000 \baselineskip=6pt
#1}\hfuzz=1pt}

\def\sss{\subsubsection}

\newcommand{\di}{\partial}

\newcommand{\ra}{\rightarrow}

\newcommand{\imply}{\Rightarrow}

\def\ssk{\smallskip}
\def\msk{\medskip}

\def\nin{\noindent}

\def\sm{\smallsetminus}

\newcommand{\diam}{\operatorname{diam}}

\newcommand{\cl}{\operatorname{cl}}
\newcommand{\inter}{\operatorname{int}}
\renewcommand{\mod}{\operatorname{mod}}

\newcommand{\tl}{\tilde}

\newcommand{\orb}{\operatorname{orb}}

\newcommand{\id}{\operatorname{id}}

\newcommand{\depth}{\operatorname{depth}}
\newcommand{\bidepth}{\operatorname{bidepth}}

\renewcommand{\d}{{\diamond}}

\newcommand{\SUBset}{\Subset}
\newcommand{\SUPset}{\Supset}

\newcommand{\eps}{{\varepsilon}}

\newcommand{\de}{{\delta}}
\newcommand{\la}{{\lambda}}
\newcommand{\La}{{\Lambda}}

\newcommand{\Om}{{\Omega}}
\newcommand{\om}{{\omega}}

\newcommand{\ba}{{\mbox{\boldmath$\alpha$} }}
\newcommand{\sba}{{\mbox{\scriptsize\boldmath$\alpha$} }}
\newcommand{\balpha}{{\mbox{\boldmath$\alpha$} }}
\newcommand{\sbb}{{\mbox{\scriptsize\boldmath$\beta$} }}
\newcommand{\bbe}{{\mbox{\boldmath$\beta$} }}

\newcommand{\sbg}{{\mbox{\scriptsize\boldmath$\gamma$} }}
\newcommand{\bgamma}{{\mbox{\boldmath$\gamma$} }}

\newcommand{\BB}{{\mathcal B}}
\newcommand{\CC}{{\mathcal C}}

\newcommand{\GG}{{\mathcal G}}

\newcommand{\KK}{{\mathcal K}}
\newcommand{\LL}{{\mathcal L}}
\newcommand{\MM}{{\mathcal M}}

\newcommand{\PP}{{\mathcal P}}

\newcommand{\RR}{{\mathcal R}}

\newcommand{\WW}{{\mathcal W}}

\newcommand{\YY}{{\mathcal Y}}

\newcommand{\A}{{\mathbb A}}
\newcommand{\C}{{\mathbb C}}

\newcommand{\D}{{\mathbb D}}

\newcommand{\N}{{\mathbb N}}

\newcommand{\R}{{\mathbb R}}
\newcommand{\T}{{\mathbb T}}

\newcommand{\U}{{\mathbb U}}

\newcommand{\Z}{{\mathbb Z}}

\newcommand{\f}{{\bf f}}
\newcommand{\g}{{\bf g}}
\newcommand{\h}{{\bf h}}
\renewcommand{\i}{{\mathbf{ i}}}

\newcommand{\n}{{\mathbf{n}}}
\newcommand{\q}{{\mathbf{q}}}

\newcommand{\s}{{\mathbf{s}}}
\newcommand{\bl}{{\boldsymbol{\lambda}}}
\renewcommand{\r}{{\mathbf{r}}}
\newcommand{\bt}{{\mathbf{t}}}
\newcommand{\psimod}{{ \mathbf{mod}\, }}
\renewcommand{\k}{{\mathbf{k}}}

\def\Bd{{\mathbf{d}}}

\def\BB{{\mathbf{B}}}

\def\Bpsi{{\boldsymbol{\Psi}}}

\def\B0{{\mathbf{0}}}
\def\BUps{{\boldsymbol{\Upsilon}}}
\def\BF{{\mathbf{F}}}
\def\BLa{{\boldsymbol{\La}}}
\def\BP{{\mathbf{P}}}
\def\BQ{{\mathbf{Q}}}

\def\BV{{\mathbf{V}}}
\def\BU{{\mathbf{U}}}
\newcommand{\BW}{{\mathbf{W}}}
\def\BY{{\mathbf{Y}}}
\def\BZ{{\mathbf{Z}}}
\def\BE{{\mathbf{E}}}

\def\Ups{{\Upsilon}}
\def\BUps{{\boldsymbol{\Upsilon}}}
\def\BLa{{\boldsymbol{\La}}}

\newcommand{\Dom}{\operatorname{Dom}}

\renewcommand{\Im}{\operatorname{Im}}

\newcommand{\from}{\mathpunct:}

\catcode`\@=12

\def\Empty{}
\newcommand\oplabel[1]{
  \def\OpArg{#1} \ifx \OpArg\Empty {} \else
  	\label{#1}
  \fi}
		
\newcommand{\comm}[1]{}
\newcommand{\comment}[1]{}

\def\makeabbrevs{%
\def\o{\omega}\def\g{\gamma}\def\G{\Gamma}\def\h{\hat}\def\d{\delta}\def\D{\Delta}%
\def\O{\Omega}\def\b{\beta}\def\l{\lambda}\def\U{\Upsilon}}

\begin{document}

\bigskip\bigskip

\title[Molecules]{A priori bounds for some\\ infinitely renormalizable quadratics:
  {\small III. Molecules.} }
\author {Jeremy Kahn and Mikhail Lyubich}
\date{\today}

\begin{abstract} 
  In this paper we prove {\it a priori bounds} for infinitely renormalizable quadratic
polynomials satisfying a ``molecule condition''. Roughly speaking, this condition ensures that
the renormalization combinatorics stay away from the satellite types. These {\it a priori bounds}
imply local connectivity of the corresponding Julia sets and the Mandelbrot set at the corresponding parameter
values.   
\end{abstract}

\setcounter{tocdepth}{1}
 
\maketitle
\thispagestyle{empty}

\def\IMSmarkvadjust{0 pt}
\def\IMSmarkhadjust{0 pt}
\def\IMSmarkhpadding{0 pt}
\def\IMSpubltext{Published in modified form:}
\def\SBIMSMark#1#2#3{
 \font\SBF=cmss10 at 10 true pt
 \font\SBI=cmssi10 at 10 true pt
 \setbox0=\hbox{\SBF \hbox to \IMSmarkhpadding{\relax}
                Stony Brook IMS Preprint \##1}
 \setbox2=\hbox to \wd0{\hfil \SBI #2}
 \setbox4=\hbox to \wd0{\hfil \SBI #3}
 \setbox6=\hbox to \wd0{\hss
             \vbox{\hsize=\wd0 \parskip=0pt \baselineskip=10 true pt
                   \copy0 \break%
                   \copy2 \break%
                   \copy4 \break}}
 \dimen0=\ht6   \advance\dimen0 by \vsize \advance\dimen0 by 8 true pt
                \advance\dimen0 by -\pagetotal
	        \advance\dimen0 by \IMSmarkvadjust
 \dimen2=\hsize \advance\dimen2 by .25 true in
	        \advance\dimen2 by \IMSmarkhadjust

%
%
  \openin2=publishd.tex
  \ifeof2\setbox0=\hbox to 0pt{}
  \else 
     \setbox0=\hbox to 3.1 true in{
                \vbox to \ht6{\hsize=3 true in \parskip=0pt  \noindent  
                {\SBI \IMSpubltext}\hfil\break
                \input publishd.tex 
                \vfill}}
  \fi
  \closein2
  \ht0=0pt \dp0=0pt
 \ht6=0pt \dp6=0pt
 \setbox8=\vbox to \dimen0{\vfill \hbox to \dimen2{\copy0 \hss \copy6}}
 \ht8=0pt \dp8=0pt \wd8=0pt
 \copy8
 \message{*** Stony Brook IMS Preprint #1, #2. #3 ***}
}

\SBIMSMark{2007/4}{December 2007}{}
\tableofcontents

\section{Introduction}

\begin{figure} \label{molecule}
  \includegraphics{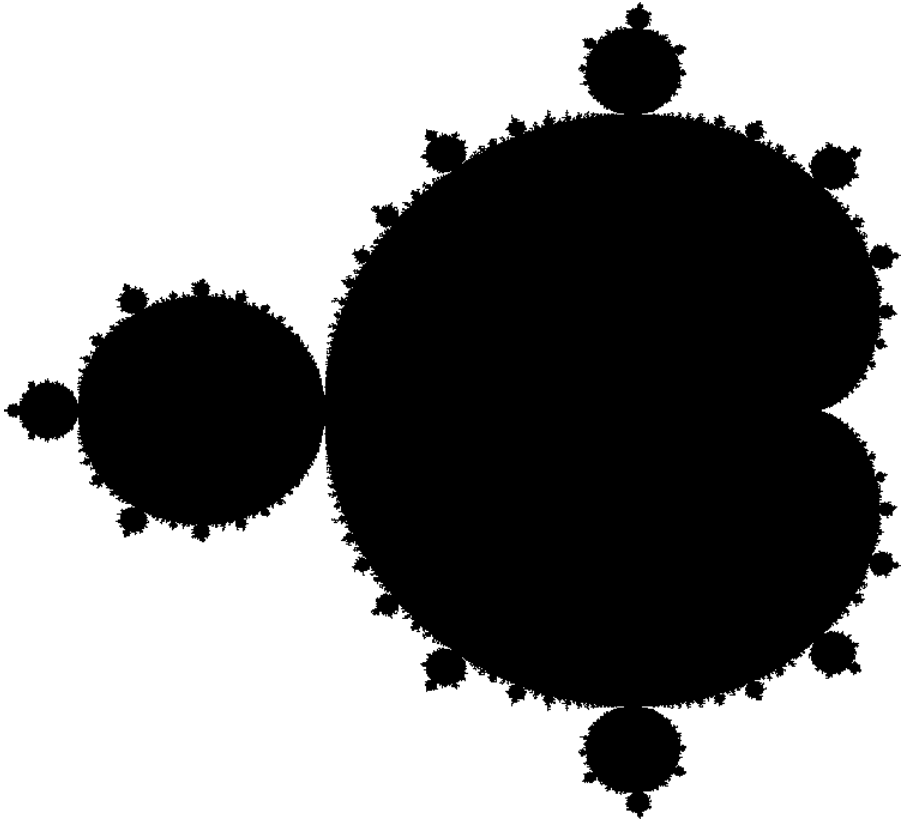}
  \caption{The central molecule of the Mandelbrot set}
\end{figure}

  The most prominent component of interior of the Mandlebrot set $M$ is the one bounded by the main cardioid. 
There are infinitely many secondary hyperbolic  components of $\inter M$ attached to it.
In turn, infinitely many  hyperbolic components are attached to each of the secondary components, etc. 
Let us take the union of all hyperbolic components of $\inter M$  obtained this way,
close it up and fill it in (i.e., add all bounded components of its complement%
\footnote{These bounded components could be only queer components of $\inter M$}).
We obtain the set called the {\it molecule} $\MM$ of $M$,
see Figure \ref{molecule}.%
\footnote{It is also called the {\it cactus}.}
In this paper we consider infinitely primitively
renormalizable quadratic polynomials satisfying a {\it molecule condition},
which means that the combinatorics of the primitive renormalization operators involved stays away from 
the molecule (see \S \ref{renorm-sec} for the precise definition in purely combinatorial terms).

An infinitely renormalizable quadratic map $f$ is said to have {\it a priori bounds}
if its renormalizations can be  represented by quadratic-like maps  $R^n f: U_n'\ra U_n$
with $\mod (U_n\sm U_n')\geq \mu >0$, $n=1,2\dots$.

Our goal is to prove the following result:

\newtheorem*{Main}{Main Theorem}.
\begin{Main}
  Infinitely renormalizable quadratic maps satisfying the molecule condition have
a priori bounds.
\end{Main}

By \cite{puzzle}, this implies:

\begin{cor}\label{JLC and MLC}  
 Let $f_c: z\mapsto z^2+c$ be an infinitely renormalizable quadratic map satisfying the molecule condition.
Then the Julia set $J(f_c)$ is locally connected, and the Mandelbrot set $M$ is locally connected at $c$.
\end{cor}

Given and $\eta>0$,
let us say that an renormalizable  quadratic map satisfies the $\eta$-{\it molecule} condition
if the combinatorics of the renormalization operators involved stays $\eta$-away from 
the molecule of $\MM$. 
 
In this paper we will deal with the case of renormalizations with sufficiently high periods.
Roughly speaking, we show that if  a quadratic-like map is nearly degenerate 
then its geometry is improving under such a renormalization.
The precise statement requires the notion of  ``pseudo-quadratic-like map'' $f$
defined in  \S \ref{pseudo-puzzle}, and its modulus, $\psimod(f)$.                   

\begin{thm}\label{high periods thm}
 Given  $\eta>0$ and $\rho\in (0,1)$, 
 there exist $\bar\mu>0$ and $\underline p\in \N$ with the following property.  
Let  $f$ be a renormalizable with period $p$ quadratic-like map
satisfying the $\eta$-molecule condition.
If $\psimod(Rf)<\bar\mu$ and $p\geq \underline p$ 
then $\psimod(f) < \rho\, \psimod (Rf)$.
\end{thm} 

The complementary case of ``bounded periods'' is dealt in \cite{K}.
 
\begin{rem}
  Theorem \ref{high periods thm} is proved in a similar way as a more special  result of \cite{decorations}.
The main difference occurs on the top level of the Yoccoz puzzle, which is modified here
so that it is associated with an appropriate periodic point rather than with the fixed point of $f$.%
\footnote{It is similar to the difference between ``non-renormalizable'' and ``not infinitely renormalizable''
cases in the Yoccoz Theorem.} 
We will focus on explaining these new elements, while only outlining the parts
that are similar to \cite{decorations}.     
\end{rem}

Let us now outline the structure of the paper.

In the next section, \S \ref{principal nest}, we lay down the combinatorial framework
for our result, the Yoccoz puzzle associated to dividing cycles, 
and formulate precisely the Molecule Condition. 

In \S \ref{pseudo-puzzle} we summarize necessary background about
{\it pseudo-quadratic-like maps} introduced in \cite{K}, 
and the {\it pseudo-puzzle} introduced in \cite{decorations}.
From now on, the usual puzzle will serve only as a combinatorial frame,
while all the geometric estimates will be made for the pseudo-puzzle.
 Only at the last moment (\S \ref{conclusion})
we return back to the standard quadratic-like context.

In \S \ref{covering lemma sec} we formulate a Transfer Prinsiple,
that will allow us to show that if $\psimod(Rf)$ is small then 
the the pseudo-modulus in between appropriate puzzle pieces is even smaller.

In \S \ref{moduli} we apply the Transfer Principle to the dynamical context.
It implies that the extremal pseuso-distance between two specific parts of the Julia set 
(obtained by removing from the Julia set the central puzzle piece $Y^1$)
is  much bigger than $\psimod(f)$ (provided the renormalization period is big). 
On the other hand, we show that under the Molecule Condition,
this pseudo-distance is comparable with $\psimod(f)$. This yields Theorem \ref{high periods thm}. 
 
\subsection{Terminology and Notation}   
  $\N =\{1,2,\dots\}$ is the set of natural numbers; $\Z_{\geq 0} = \N \cup \{0\}$;
$\D=\{z:\,  |z|<1\}$ is the unit disk, and $\T$ is the unit circle.\\

    A  {\it topological disk} means a simply connected  domain in $\C$.
A {\it continuum} $K$ is a connected closed subset in $\C$. 
It is called {\it full} if all components of $\C\sm K$ are unbounded.   

For subsets $K, Y$ of a topological space $X$,   
notaion $K\Subset Y$ will mean (in a slightly non-standard way) 
that the closure  of $K$ is a compact contained in $\inter Y$.   

  We let $\orb(z)\equiv \orb_g(z)= (g^n z)_{n=0}^\infty$ be the {\it orbit} of $z$ under a map $g$.

  Given a map $g: U\ra V$ and an open topological disk $D\subset V$, 
components of $g^{-1}(D)$ are called {\it pullbacks}  of $D$ under $g$. 
If the disk $D$ is closed, we define pullbacks of $D$ as the closures of the pullbacks of $\inter D$.%
\footnote{Note that the pullbacks of a closed disk $D$ can touch one another, 
so they are not necessarily connected components of $g^{-1}(D)$.} 
In either case, 
given a connected set $X\subset g^{-1}(\inter D)$, we let $g^{-1} (D)|X$ be the pullback of $D$
containing $X$.    


\subsection{Acknowledgement}
We thank Scott Sutherland for help with making Figure \ref{molecule}.
This work has been partially supported by the NSF and NSERC.

\section{Dividing cycles, Yoccoz puzzle, and renormalization}\label{principal nest}

Let $(f: U'\ra U)$ be a quadratic-like map.
We assume that the domains $U'$ and $U$ are smooth disks,
$f$ is even,  
and we  normalize $f$ so that $0$ is its critical point.  
We let $U^m = f^{-m}(U)$. The boundary of $U^m$ is called  
the {\it equipotential of depth $m$}.   

By means of straightening, we can define external rays for $f$. They form a foliation of $U\sm K(f)$
transversal   to the equipotential $\di U$. Each ray is labeled by its {\it external angle}. 
These rays will play purely combinatorial role, so particular choice of the straightening
is not important. 
  
\subsection{Dividing cycles and associated Yoccoz puzzles}\label{sec: first king} 

  Let us consider a repelling periodic point $\gamma$ of period $\bt$ and the correspoding cycle 
$\bgamma=\{f^k\gamma\}_{k=0}^\bt$. This point (and the cycle) is called {\it dividing}
if there exist at least two rays landing at it.  
For instance, the landing point of the zero ray is a non-dividing fixed point, 
while the other fixed point is dividing (if repelling). 

{\it In what follows, we assume that $\gamma$ is dividing.}
Let $\RR(\gamma)$ (resp., $\RR(\bgamma)$) stand for the family of rays landing at $\gamma$ (resp., $\bgamma$).
Let ${\bf s} =\# \RR(\gamma)$ and let $\r= \bt\s= \# \RR(\bgamma)$. 
These rays divide $U$ into $\bt(\s-1)+1$ closed topological disks $Y^0(j)\equiv Y^0_\sbg(j)$ called 
{\it Yoccoz puzzle pieces of depth $0$}.

{\it Yoccoz puzzle pieces $Y^m(j)\equiv Y^m_\sbg(j)$ of depth $m$} are  defined as the 
pullbacks of $Y^0(i)$ under $f^m$. 
They tile the neighborhood of $K(f)$ bounded by the equipotential $\di U^m $.
Each of them is bounded by finitely many arcs of this equipotential and finitely many  
external rays of $f^{-m}(\RR(\bgamma))$. 
We will also use notation $ Y^m(z)$ for the puzzle piece $Y^m(j)$ containing $z$ in its interior. 
If $f^m(0)\not\in \bgamma$, then there is a well defined {\it critical} puzzle piece  $Y^m\equiv Y^m(0)$.
The critical puzzle pieces are nested around the origin:
$$
   Y^0\supset Y^1\supset Y^2 \dots \ni 0. 
$$
Notice that all $Y^m$, $m\geq 1$, are symmetric with respect to the origin.

Let us take a closer look at some puzzle piece $Y=Y^m(i)$.
Different arcs of $\di Y$  meet at the {\it corners} of $Y$.
The corners where two external rays meet will be called {\it vertices} of $Y$;
they are $f^m$-preimages of $\bgamma$.
Let $K_Y= K(f)\cap Y$. It is a closed connected set that meets the boundary $\di Y$
at its  vertices. Moreover, the external rays meeting at a vertex $v\in \di Y$ chop off from
$K(f)$ a continuum $S_Y^v$, the component of $K(f)\sm \inter Y$ containing $v$. 

Let $\YY_\sbg$ stand for the family of all puzzle pieces $Y_\sbg^m(j)$. 

Let us finish with an obvious observation that will be constantly exploited:

\begin{lem}\label{Subset}
  If a puzzle piece $Y^n(z)$ of $\YY_\bgamma$ does not touch the cycle $\bgamma$
then $Y^n(z)\Subset Y^0(z)$.
\end{lem}

\subsection{Renormalization associated with a dividing cycle}\label{renorm-sec}

\begin{lem}[see \cite{Th,M-rays}]\label{cv pp}
  The puzzle piece $X^0 \equiv Y^0(f(0))$ of $\YY_\bgamma$ 
containing the critical value has only one vertex, 
and thus is bounded by only two external rays (and one equipotential).
\end{lem} 
  
In what follows, 
$\gamma$ will denote the point of the cycle $\bgamma$ such that $f(\gamma)$ is the vertex of $X^0$.
Notice that $f(0)\in \inter X^0$ for otherwise $f(0)=f(\gamma)$, which is impossible since $0$ is the
only preimage of $f(0)$. 

Since the critical puzzle piece $Y^1$ is the pullback of $X^0$ under $f$, it has two vertices, 
$\gamma$ and $\gamma'= -\gamma$, and is bounded by four rays,
two of them landing at $\gamma$ and two landing at $\gamma'$.%
\footnote{We will usually say that ``a puzzle piece is bounded by several external rays'' 
without mentioning equipotentials that also form part of its boundary.} 

\begin{lem}(see \cite{D})
  Let $X^\r \equiv Y^\r (j)\subset X^0$  be the puzzle piece attached to the boundary of $X^0$.
Then $f^\r : X^\r \ra X^0$ is a double branched covering.
\end{lem}
 
\begin{pf}
    Let $\CC^0$ be the union of the two rays that bound $X^0$,
and let $\CC^\r  $ be $\CC^0$ cut by the equipotential $\di U^\r$.
Let us orient $\CC^0$ and then induce the orientation to $\CC^\r $.
    Since $f^\r : \CC^\r \ra \CC^0$ is an orientation preserving homeomorphism, 
it maps $X^\r $ onto $X^0$. 

Since for $m=1,\dots, \r - 1$,  the arcs $f^m(\CC^\r)\subset \di (f^m X^\r )$ are disjoint  from
$\inter X^0$, the puzzle pieces $f^m X^\r $ are not contained in $X^0$.
Since they have a bigger depth than $X^0$, they are disjoint from $\inter X^0$.
It follows that all the puzzle pieces $f^m (X^\r )$, $m=0,1,\dots, \r -1$, 
have pairwise disjoint interiors. 
(Otherwise $f^m(X^\r)\supset f^n (X^\r)$ for some $ \r> m > n\geq 0$,
and applying $f^{\r-m}$, we would conclude that $X^0\supset f^{\r-n+m} (X^\r)$.)    

Moreover, the puzzle piece  $f^{\r-1}X^\r=Y^1$ is  critical
since $Y^1$ is the only pullback of $X^0$ under $f$. 
Hence the puzzle pieces $f^m(X^\r)$, $m=0,1,\dots, \r-2$, do not contain $0$. 
It follows that  $\deg(f^\r : X^\r \ra X^0)=2$.
\end{pf}

\begin{figure}[htbp]
\begin{center}
{\makeabbrevs
\input{figures/pullback2.pstex_t}
\caption{}
\label{pullback}}
\end{center}
\end{figure}

\begin{cor}\label{renorm}
  If $f(0)\in \inter X^\r $ then the puzzle piece $Y^{\r+1}$ has four vertices,
and the map $f^\r : Y^{\r+1}\ra Y^1$ is a double branched covering.  
\end{cor}

   Let $\Theta(\bgamma)\subset \T$ be the set of external angles of the rays of $\RR(\bgamma)$.
There is a natural equivalence relation on $\Theta(\bgamma)\subset \T$:
two angles are equivalent if the corresponding rays land at the same periodic point. 
Let us consider the hyperbolic convex hulls of these equivalence classes 
(in the disk $\D$ viewed as the hyperbolic plane). 
The union of the boundaries of these convex hulls 
is a  finite lamination $\PP=\PP(\bgamma)$ in $\D$ which is also called the 
{\it periodic ray portrait}. One can characterize all possible ray portraits that appear in this way
(see \cite{M-rays}).

\begin{definition} 
A map $f$ is called $\PP$-{\it renormalizable} 
(or, ``renormalizable with combinatorics $\PP$'')
if $f(0)\in \inter X^\r$ and $f^{\r m}(0)\in  Y^{\r+1}$  
for all $m=0,1,2,\dots$. 
In this case, the double covering $f^\r : Y^{\r+1}\ra Y^1$ is called
the {\it renormalization} $R_\PP f = R_\sbg f$ of $f$ (associated with the cycle $\bgamma$).
The corresponding {\it  little (filled) Julia set} $\KK = K(R f)$ 
is defined as 
$$
        \{ z: f^{\r m}z\in Y^{\r + 1} ,\quad  m=0,1,2\dots\}
$$
If the little Julia sets $f^m \KK$, $k=0,1\dots, \r-1$, 
are pairwise disjoint, then the renormalization is called {\it primitive};
otherwise it is called {\it satellite}. 

In case $\bgamma$ is the dividing fixed point of $f$,
the map  is also called {\it immediately renormalizable}.
(This is a particular case of the satellite renormalization.)
\end{definition}

\begin{rem} The above  definition of renormalization is not quite standard since the map 
$R_\sbg f$ is not quadratic-like. To obtain the usual notion of renormalization, 
one should thicken the domain of $R_\sbg f$  a bit to make it quadratic-like (see \cite{D,M}). 
This thickenning does not change the Julia set, so $K(Rf)$ possesses all the properties of 
quadratic-like Julia sets. In particular, it has two fixed points, one of which is either
non-repelling or dividing.

Note also that in the case when $f^r(0)=\gamma'$ the puzzle piece
$Y^{\r+1}$ degenerates (is pinched at $0$), but this does not effect any of 
further considerations. 
\end{rem}

Given a periodic ray portrait $\PP$, 
the set of parameters $c\in \C$ for which the quadratic polynomial $P_c$ is $\PP$-renormalizable form a 
{\it little copy} $M_\PP$ of the Mandelbrot set (``$M$-copy''). Thus, there is one-to-one correspondence
between the admissible ray portraits and the little M-copies.
So, one can encode the combinatorics of the renormalization by the little $M$-copies themselves. 

\subsection{Molecule Condition}

The  molecule $\MM$ defined in the Introduction consists of the quadratic maps which are:

\nin\ssk $\bullet$
 either finitely many times renormalizable, all these renormalizations are satellite,
  and the last renormalization has a non-repelling cycle;

\nin\ssk $\bullet$
  or  infinitely many times renormalizable, with all the renormalizations satellite.

\ssk The molecule condition that we are about to introduce will  ensure that our map $f$
  has frequent ``qualified'' primitive renormalizations.
Though  $f$ is allowed to be satellite renormalizable once in a while,
{\it we will   record only the primitive renormalizations}.
They  are naturlly ordered according to their periods,
$1=p_0 <p_1<p_2<\dots$, where  $p_i$ is a multiple of $p_{i-1}$. 

Along with these ``absolute'' periods of the primitive renormalizations,
we will consider {\it relative periods} $\tl p_i = p_i/p_{i-1}$ and the corresponding $M$-copies
$\tl M_i$ that encode the combinatorics of $R^i f$ as the renormalization of $R^{i-1} f$. 


Given an $\eta>0$, we say that a sequence of primitively renormalizable quadratic-like maps $f_i$
satisfies the $\eta$-{\it molecule condition}
if the corresponding  $M$-copies $M_i$  stay $\eta$-away from the molecule $\MM$ 
(the latter is defined in the Introduction).
We say that $\{f_i\}$ satisfies the {\it molecule condition} if it does it for some $\eta>0$.
    
An infinitely primitively  renormalizable map $f$  satisfies the $\eta$-molecule condition
if the sequence of its primitive renormalizations $R^i f$ does
(i.e.,  the corresponding relative copies $\tl M_i$ stay $\eta$-away from the molecule $\MM$).   
(And similarly, for  the non-quantified molecule condition.)

\msk There is, however, a more specific combinatorial way to describe the molecule condition. 

Let us consider a quadratic-like map $f$ with straightening $P_c$, $c\in M$.
We are going to associate  to  $f$ (in some combinatorial region, and with some choice involved)
three combinatorial parameters, $(\r,\q,\n)$ (``period, valence, and escaping time'')
whose boundedness will be equivalent to the molecule condition.
 
\ssk
{\it Assume first that $f$ admits a dividing cycle $\bgamma$ with the ray portrait $\PP$
with $\r$ rays}. 
This happens if and only if $c$ belongs to the parabolic limb of $M$ cut off by
two external rays landing at an appropriate parabolic point. 

\ssk
On the central domain of $\C\sm (\RR(\bgamma)\cup \RR(\gamma'))$,
$f^\r$ has a unique fixed poin $\alpha$. 
{\it Next, we assume that $\alpha$ is repelling and there are $\q$ rays landing at it.} 

\ssk {\it Assume next that the finite orbit $f^{\r j(0)}$, $j=1,\dots, \q\n-1$,
      does not escape the central domain of  $\C\sm (\RR(\gamma)\cup \RR(\gamma'))$}.
       This happens if and only if 
  $c$ lies outside certain {\it decorations} (see \cite{decorations})
  of the above parabolic limb. In particular, this happens if $f$ is
  {\it satellite  $\PP(\gamma$)-renormalizable}.   

\ssk {\it Finally, assume that  $\n$ is the first moment $n$ such that  $f^{\r\q n}(0)$ escapes the
  central domain of $\C\sm (\RR(\alpha)\cup \RR(\alpha'))$}. 
   This happens if and only if 
  $c$ belongs to the union  of $2^\n$ {\it decorations} 
  inside the above parabolic limb. In particular, the map $f$ 
  {\it is not  $\PP(\balpha$)-renormalizable}.   

\ssk 
 Under the above assumptions, 
  we say that $f$ satisfies the $(\bar\r, \bar\q, \bar\n)$-{\it molecule condition} if
  there is a choice of $(\r, \q,\n)$ with $\r\leq \bar\r$, $\q\leq \bar \q$ and $\n\leq \bar \n$.

\begin{lem}
     The $(\r,\q,\n)$-molecule condition is equivalent to the $\eta$-molecule condition.%
\footnote{In the sense that if $f$ satisfies $(\r,\q,\n)$-molecule condition then it satisfies
 $\eta$-molecule condition with some $\eta=\eta(\r,\q,\n)$, and the other way around.}  
\end{lem}

\begin{pf}
  If $f$ satisfies $(\r,\q,\n)$-molecule condition then $c$ belongs to the finite union of decorations.
Each of them does not intersect the molecule $\MM$, so $c$ stays some distance $\eta$ away from $\MM$. 

Vice versa, assume there is a sequence of maps $f_i$ satisfying the $\eta$-molecule condition,
but with  $(\r,\q,\n)\to \infty$ for any choice of $(\r,\q,\n)$.
Let us select a convergent subsequence $c_i\to c$. Since $c\not\in \MM$, 
there can be only finitely many hyperbolic components $H_0,H_1,\dots,H_m$ of $\MM$
such that $H_0$ is bounded by the main cardioid, $H_{k+1}$ bifurcates from $H_k$,
and $f$ is $m$ times immediately  renormalizable with the corresponding combinatorics.   
Let us consider the two rays landing at the last bifurcation point (where $H_m$ is attached to $H_{m-1}$),
and the corresponding parabolic limb of the Mandelbrot set.%
\footnote{If $m=0$ then we consider the whole Mandelbrot set. } 
This parabolic point has certain period $\r$. 

Since the quadratic polynomial $P_c$ is not immediately renormalizable any more,
the corresponding cycle $\balpha$ is repelling with $\q$ rays landing at each of  its periodic points,
and there is some escaping time $\n$. 
So, $P_c$ satisfies $(\r,\q,\n)$-molecule condition.
Since this condition is stable under perturbations, $P_{c_i}$ satisfy it as well -- contradiction.  
\end{pf}


{\it In what follows we assume that parameters $\r,\q,\n$ are well defined for a map $f$ under consideration,}
so in particular, we have two dividing cycles, $\bgamma$ and $\balpha$. 
We let $\k=\r\q\n$.
Let us state for the record the following well-known combinatorial property:

\begin{lem}\label{zeta}
   The point $\zeta=f^\k(0)$ is separated from $\alpha$ and $0$ by the rays landing at $\alpha'$.%
\footnote{while the latter two points are not separated.}
\end{lem}

 \subsection{Combinatorial separation between $\gamma$ and $\balpha$}\label{comb sep} 
  

Along with the puzzle $\YY_\sbg$ associated with $\bgamma$, let us consider the puzzle $\YY_\sba$
associated with $\balpha$. 
The critical puzzle piece $Y_\sba^1$ has two vertices, $\alpha$ and $\alpha'$,
and is bounded by four external rays landing at these vertices. Let $\CC$ 
be the union of the two external rays of $\di Y^1_\sba$ landing at $\alpha$, 
and let $\CC'$  be the symmetric pair of external rays landing at $\alpha'$. 
     
Recall that $\bt$ stands for the period of $\bgamma$ and $\s$ stands for $\# \RR(\bgamma)$. 

\begin{lem}\label{separation}
  There exist inverse branches $f^{-\bt m} |\, \CC'$, $m=0,1,\dots, \s-1$, 
such that the union of the arcs $f^{-\bt m} (\CC')$  separates $\gamma$ from 
the cycle $\balpha$ and the co-cycle $\balpha'$
(except that $\alpha'\in \CC'$). 
\end{lem}

\begin{pf}
Let us pull the puzzle piece $Y_\sbg^1$ along the orbit $\bgamma$ 
(or equivalently, along the orbit  $\balpha$). 
By Corollary \ref{renorm}, 
the corresponding  inverse branches 
$$
   f^{-m} : Y^1_\sbg \ra f^{-m}(Y^1_\sbg) = f^{\r-m} Y^{\r+1}_\sbg,\quad  m=0,1,\dots, \r-1,
$$ 
have disjoint interiors.  
Hence each of these puzzle pieces contains exactly one point of $\balpha$. 
Moreover, non of these puzzle pieces except $Y^1_\sbg$ may intersect $\balpha'$
(for otherwise, its image would contain two points of $\balpha$).   

It follows from standard properties of quadratic maps 
that the arc $\CC'$ {\it separates} 
$\gamma$ (which is the non-dividing fixed point of $R_\sbg f$) from $\gamma'$ and
$\alpha$ (which is the dividing fixed point of $R_\sbg f$). 
Hence the arcs $f^{-\bt m} (\CC')$, 
separate $\gamma$ from  $f^{-\bt m}(\gamma')$ and $f^{-\bt m}(\alpha)$. 

Since each of the puzzle pieces $f^{-\bt m} ( Y^1_\sbg)$,  $m=0,1,\dots, \s-1$, has two vertices
($\gamma$ and $f^{-\bt m}(\gamma')$) and  their union forms a neighborhood of $\gamma$,
the rest of the Julia set is separated from $\gamma$ by the union of arcs
$f^{-\bt m}(\CC')$, $m=0,1,\dots, \s-1$. 
It follows that this union separates $\gamma$ from the whole cycle $\balpha$, 
and from the co-cycle $\balpha'$. 
\end{pf}

Together with Lemma \ref{Subset} this yields: 

\begin{cor} We have:
  $Y_\sba^{\r}(\gamma) \SUBset Y^0_\sba (\gamma)$ 
\end{cor}

\begin{pf}
   The puzzle piece $Y_\sba^{\r}(\gamma)$ contains $\gamma$ and does not cross the arcs 
$f^{-\bt m} (\CC')$, $m=0,1,\dots, \s-1$, from Lemma \ref{separation}.
Hence it does not intersect $\balpha$, and the the conclusion follows from Lemma \ref{Subset}.
\end{pf}

 \subsection{A non-degenerate annulus} 

 {\it In what follows we will be dealing only with the puzzle $\YY_\ba$, so we will skip the label $\ba$
in notation.} 
  
Since by Lemma \ref{zeta}, the point $\zeta =f^\k (0)$ is separated from $\alpha$ by $\CC'$,
the union of arcs $f^{-\bt m} (\CC')$, $m=0,1,\dots, \s-1$, from Lemma \ref{separation}
separates $\zeta$ from the whole cycle $\balpha$.
By Lemma \ref{Subset},  $Y^{\r}(\zeta)\Subset Y^0(\zeta)$.
Pulling this back by $ f^\k$,  we conclude:

\begin{lem}
We have: 
       $Y^{\k+\r} \Subset  Y^\k$.
\end{lem}

We let $E^0= Y^{\k+\r}$. 

\subsection{Buffers attached to the vertices of $P=Y^{\r\q (\n-1)+1}$}\label{P-sec}

 Let us consider a nest of critical puzzle pieces
$$
   Y^1 \supset  Y^{ \r\q +1 } \supset Y^{2\r\q+1 } \supset\dots \supset Y^{\r\q (\n-1)+1} = P. 
$$
 Since $f^{\r\q} : Y^{\r\q+1}\ra Y^1 $ is a double branched covering such that
$$
   f^{\r\q m}(0)\in Y^1, \quad m=0,1,\dots, \n-1,
$$
the puzzle piece $Y^{\r\q k+1}$  is mapped by 
 $f^{\r\q}$ onto $Y^{\r\q (k-1)+1}$ as a double branched covering, $k=1\dots,\n-1 $. 
However, since $f^\k(0)=f^{\r\q\n}(0)\not\in Y^1$, there are two non-critical puzzle pieces of depth 
$\k+1$ mapped univalently onto $P$ under $g^{\q}$.  
One of these puzzle pieces, called $Q_L$, is attached to the point $\alpha$,
another one, called $Q_R$, is attached to $\alpha'$.
The following Lemma is similar to Lemma  2.1 of \cite{decorations}:

\begin{lem}\label{Qv}
 For any vertex of $P$, there exists a puzzle piece $Q^v\subset P$ of
 depth $\r(2\n-1)\q + 1$ attached to the boundary rays of $P$ landing at $v$
 which is a univalent $f^\k$-pullback of $P$. Moreover, these puzzle pieces are
 pairwise disjoint.
\end{lem}

\subsection{Modified principle nest}\label{modified nest}
Until now, the combinatorics of the puzzle depended only on the parameters $(\r,\q,\n)$.
Now we will dive  into the deeper waters.

Let $l$ be the first return time of $0$ to $\inter E^0$ and let $E^1= Y^{\k+ \r +l}$
be the pullback of $E^0$ along the orbit $\{f^m (0)\}_{m=0}^l$.
Then $f^l: E^1\ra E^0$ is a double branched covering.

\begin{cor}\label{non-degenerate annulus}
     We have: $E^1\SUBset  E^0$.
\end{cor}

\begin{pf}
  Since $\{f^m(0)\}_{m=1}^{\r-1}$ is disjoint from $Y_\sbg^1\supset Y_\sba^1 \supset E^0$,
we have: $l\geq \r$. Hence 
$$
      f^l(E^0)\supset f^\r(E^0) \Supset E^0 = f^l(E^1),
$$
and the conclusion follows. 
\end{pf}

Given two critical puzzle pieces $E^1\subset \inter E^0$, we can construct 
the {\it (Modified) Principle Nest} of critical puzzle pieces 
$$
  E^0\SUPset E^1\SUPset E^1\SUPset\dots\SUPset E^{\chi-1}\SUPset E^\chi 
$$
as described in  in \cite{KL-high}.                        
It comes together with quadratic-like maps $g_n: E^n\ra E^{n-1}$. 

If the map $f$ is renormalizable then  the Principle Nest terminates at some level $\chi$. 
In this case, the last quadratic-like map $g_\chi: E^\chi\ra E^{\chi-1}$  has connected Julia set
that coincides with  the Julia set of the  renormalization $R_\sbb f$, where
$\bbe$ is the $f$-orbit of the non-dividing fixed point $\beta$ of $g_\chi$.  
The renormalization level $\chi$ is also called the {\it height} of the nest.


\subsection{Stars}\label{stars}


Given a vertex $v$ of some puzzle piece of depts $n$,
let $S^n(v)$ stand for the union of the puzzle pieces of depth $n$ attached to $v$
(the ``star'' of $v$).
Given a finite set ${\bf v}=\{v_j\}$ of vertices $v_j$, we let
$$
       S^n({\bf v})=\bigcup_j S^n(v_j). 
$$

  Let us begin with an obvious observation that follows from Lemma~\ref{Subset}:

\begin{lem}\label{compactly conained}
  If a puzzle piece $Y^n(z)$ is not contained in $S^n(\balpha)$
then $Y^n(z)\Subset Y^0(z)$.
\end{lem}

\begin{lem}\label{single point}
  For $\bl=\k + 1+2\r$, the stars $S^{\bl}(\alpha_j)$ do not overlap and do not contain the critical point. 
\end{lem}

\begin{pf}
Let us consider the  curves $\CC$ and $f^{-\bt m} \CC'$, $m=1,\dots, \s-1$,  from Lemma \ref{separation}. 
They separate $\alpha'$ from all points of $\balpha\cup \balpha'\sm \{\alpha\}$. 
Furthermore, since $f^{\k}(0)$ is separated from $0$ by $\CC'$,
there is a lift $\Gamma$ of $\CC'$ under $f^{\k}$ that separates $\alpha'$ from $0$ 
and hence from $\alpha$.  It follows that the curves $\Gamma$ and $f^{-\bt m} \CC'$, $m=1,\dots, \s-1$,
separate $\alpha'$ from all points of $\balpha\cup \balpha'$. 
Since the maximal depth of these curves is  $\bl=\k +1$ (which is the depth of $\Gamma$),
the star $S^{\k+1}(\alpha')$ does not overlap with the interior of the stars $S^{\k+1}(a)$
for all  other $a\in \balpha\cup\balpha'$.

By symmetry, the same is true for the star $S^{\k+1}(\alpha)$.
Since these stars  do not contain $0$,
the pullback of $S^{\k+1+\r}(\alpha)$  under $f^\r$ (along $\balpha$) is compactly contained in 
its interior, $\inter S^{\k+1+\r}(\alpha)$.
It follows that $S^{\k+1+\r}(\alpha)$ does overlap with the stars
$S^{\k+1}(a)$ for all  other $a\in \balpha\cup\balpha'$.

Pulling this star once more around $\balpha$, 
we obtain a disjoint family of stars.
Hence all the stars $S^{\k+1+2\r}(a)$,  $a\in \balpha\cup\balpha'$,
are pairwise disjoint. 
\end{pf}


\subsection{Geometric puzzle pieces}

In what follows  we will deal with more general puzzle pieces.

Given  a puzzle piece $Y$, of depth $m$, 
let $Y[l]$ stand for a Jordan disk bounded by the same external rays
as $Y$ and arcs of equipotentials of level $l$ (so $Y[m]=Y$).
Such a  disk will be called a puzzle piece of {\it bidepth} $(m,l)$.
  
A {\it geometric puzzle piece} of bidepth $(m,l)$ is a closed Jordan domain which is 
the union of several puzzle pieces of the same bidepth. 
As for ordinary pieces, a pullback of a geometric puzzle piece of bidepth $(m,l)$ 
under some iterate $f^k$ is a geometric puzzle piece of bidepth $(m+k, l+k)$.
Note also that if $P$ and $P'$ are geometric puzzle pieces with%
\footnote{the inequality between bidepths is understood componentwise}
$\bidepth P \geq \bidepth P'$ 
and $K_P\subset K_{P'}$ then $P\subset P'$.

The family of geometric puzzle pieces of bidepth $(m,l)$ will be called $\YY^m[l]$.

Stars give examples of geometric puzzle pieces. 
Note that pullback of a star is a geometric puzzle piece as well 
but it is a star only if the pullback is univalent.

\comm{*******************8
Any puzzle piece $Y\in \YY^m[l]$ admits the following combinatorial representation.
Let $\theta_i$ be the cyclically ordered  angles of the external rays $\RR_i$ that bound $Y$. 
Let us consider the straight rays $R_i$ in $\C\sm \D$ of angles $\theta_i$
truncated by the circle $\T_r$ of radius $r=1/2^l$. 
If two consecutive rays, $\RR_i$ and $\RR_{i+1}$, land at the same vertex of $Y$,
let is connect $R_i$ to $R_{i+1}$ with a hyperbolic geodesic in $\D$. 
Otherwise $\RR_i$ and $\RR_{i+1}$ are connected with an equipotential arc. 
Then let us connect $R_i$ to $R_{i+1}$ with the appropriate arc of $\T_r$.
We obtain a Jordan curve that bounds the {\it combinatorial model} $M_Y$ of $Y$.

 The arcs $\om_i$ of  $\T\cap M_Y$ correspond to  the 
``external arcs'' of the Julia piece $K_Y$.  
 They have length $2\pi \la$, where $\la$ is called the {\it combinatorial length}
of the corresponding external arc of $K_Y$. 
For a geometric puzzle piece $Y$ of depth $m$,
the combinatorial length of its external arcs is at least $c(\PP)2^{-m}$,
where $c(\PP)$ is the minimum of the combinatorial lengths of the external arcs of the
puzzle pieces of depth 0. 
******************************************}

\msk

\begin{lem}\label{qn-pullbacks}
  Given a point $z\in \inter S^1(\alpha)\cup \inter S^1(\alpha')$ such that $f^{\k} z \in\inter S^1(\alpha)$,
let $P = f^{-\k}(S^1(\alpha))|z$.  Then $P\subset S^1(\alpha)$ or $P\subset S^1(\alpha')$.
\end{lem}

\begin{pf}
  Notice that $\alpha$ and $\alpha'$ are the only points of $\balpha\cup \balpha'$ 
contained in $\inter S^1(\alpha)\cup \inter S^1(\alpha')$.
Since $f^{\k}$ maps $\tl\balpha=\balpha\cup \balpha'\sm \{\alpha, \alpha'\}$
to $\balpha\sm \{\alpha\}$,
no point of $\tl\balpha$ is contained in $\inter P$.
Hence $P$ is contained in $ S^1(\alpha)\cup  S^1(\alpha')$.

But by construction of $Y^{\k+1}$, the interior of $f^\k (Y^{\k+1})$ does not overlap with $S^1(\alpha)$.
Hence $\inter Y^{\k+1}$ does not overlap with $P$. 
But $S^1(\alpha)\cup S^1(\alpha')\sm  \inter Y^{\k +1}$ 
consists of two components, one inside $S^1(\alpha)$ and the other inside $S^1(\alpha')$. 
Since $P$ is connected, it is contained in one of them.
\end{pf}

\section{Pseudo-quadratic-like maps and pseudo-puzzle}\label{pseudo-puzzle}

In this section, we will summarize the needed background on pseudo-quadratic-like maps
and pseudo-puzzle. The details can be found in \cite{K,decorations}. 

\subsection{Pseudo-quadratic-like maps}

Suppose that $\BU'$, $\BU$ are disks,
$i:  \BU' \to \BU$ is a holomorphic immersion,
and $f: \BU' \to \BU$ is a degree $d$ holomorphic branched cover.
Suppose further that there exist full continua $K \Subset \BU$ and $K'\Subset \BU'$
such that $K' = i^{-1}(K) = f^{-1}(K)$.
Then we say that $F=(i,f):  \BU'\ra (\BU, \BU)$  is a  
$\psi$-{\it quadratic-like}
 ($\psi$-{\it ql}) map  with filled Julia set $K$.
We let
$$
   \psimod(F)= \psimod(f) = \mod(\BU\sm K). 
$$ 

\begin{lem}\label{ql extension}
Let $F=(i,f)\from\BU'\ra \BU$ be a $\psi$-ql map of degree $d$ 
 with filled Julia set $K$.
Then $i$ is an embedding in  a neighborhood of $K' \equiv f^{-1}(K)$, 
and the map  $g\equiv  f\circ i^{-1}\from U'\ra U$ near $K$ is quadratic-like.

\ssk Moreover, the domains $U$ and $U'$ can be selected in such a way that
$$ \mod(U\sm i(U'))\geq \mu (\psimod(F) >0.$$ 
\end{lem}

There is a natural $\psi$-ql map $\BU^n\ra \BU^{n-1}$, 
the ``restriction'' of $(i,f)$ to $\BU^n$. 
Somewhat loosely, we will use the same notation $F=(i,f)$ for this restriction.  

 Let us normalize  the $\psi$-quadratic-like  maps under consideration 
so that $\diam K'= \diam K=1$, both $K$ and $K'$ contain 0 and 1,  
$0$ is the critical point of $f$, and $i(0)=0$.
Let us endow the space of $\psi$-quadratic-like maps (considered up to independent rescalings in the domain and the range)
with the Carath\'eodory topology.
In this topology, a sequence of normalized maps $(i_n, f_n): \BU_n'\ra\BU_n$ converges to $(i,f): \BU'\ra \BU$ if
the pointed domains $(\BU_n', 0)$ and $(\BU_n, 0)$ converge to $\BU'$ and $\BU$ respectively,
and the maps $i_n$, $f_n$ converge respectively to $i$, $f$, 
uniformly on compact subsets of  $\BU'$.

\begin{lem}\label{compactness}
  Let $\mu>0$. Then the space of $\psi$-PL maps $F$ with connected  Julia set
and  $\psimod(F)\geq \mu$ is compact.
\end{lem}  

To simplify notation, we will often refer to $f$ as a ``$\psi$-ql map'' 
keeping $i$ in mind implicitly.

\subsection{Pseudo-puzzle}\label{psi-puzzle}

\subsubsection{Definitions}
   Let $(i,f) : \BU'\ra \BU$ be a $\psi$-ql map.
By Lemma~\ref{ql extension}, it admits a quadratic-like restriction $U'\ra U$ to a neighborhood of
its (filled) Julia set  $K=K_\BU$.
Here $U'$ is embedded into $U$, so we can identify $U'$ with $i(U')$
and $ f: U'\ra U $ with  $f \circ i^{-1} $. 

Assume that $K$ is connected and both fixed points of $f$ are repelling.
Then we can cut $U$ by external rays landing at the $\alpha$-fixed point and consider the corresponding Yoccoz puzzle.

Given a (geometric) puzzle piece $Y$ of bidepth $(m,l)$,
recall that $K_Y$ stands for $Y\cap K(f)$. 
Let us consider the topological annulus $A=\BU^l  \sm K(f)$ and its universal covering $\hat A$. 
Let $Y_i$ be the components of $Y\sm K_Y$. There are finitely many of them,
and each $Y_i$ is simply connected. Hence they can be embedded into $\hat A$.
Select such an embedding $e_i: Y_i\ra \hat A_i$ where $\hat A_i$ stands for a copy of $\hat A$. 
Then glue the $A_i$ to $Y$ by means of $e_i$, i.e., let $\BY= Y \sqcup_{e_i} \hat A_i.$
This is  the {\it pseudo-piece} (``$\psi$-piece'') associated with $Y$.
Noter that the Julia piece $K_Y$ naturally embeds into $\BY$.

\begin{lem}\label{psi-maps}
\begin{itemize}
\item [(i)]  Consider two puzzle pieces $Y$ and $Z$ such that the map $f: Y\ra Z$ is a branched covering
of degree $k$ (where $k=1$ or $k=2$ depending on whether $Y$ is off-critical or not).
Then there exists an induced map $\f: \BY \ra \BZ$
which  is a branched covering of the same degree $k$.

\item [(ii)] Given two puzzle pieces $Y\subset Z$, the inclusion $i: Y\ra Z$ extends to an immersion $\i: \BY\ra \BZ$.
\end{itemize}
\end{lem}

\comm{******************************************
\sss{Moduli}
 Given two puzzle pieces $Z\Subset Y$, we let 
$$
  \psimod (Y,Z) = \mod(\BY\sm K_Z).
$$

Lemma \ref{psi-maps} implies:

\begin{lem}\label{psi-mod}
\begin{itemize}
  \item [(i)]  Consider two pairs of puzzle pieces $(Y',Y)$ and $(Z',Z)$ 
such that the map $f: (Y',Y)\ra (Z',Z)$ is a branched covering
of degree $k$ (on both domains).
Then $$\psimod(Z',Z)=k\, \psimod(Y',Y).$$

\item [(ii)] Given a nest of three puzzle pieces $W\subset Z\subset  Y$, 
we have $$\psimod(Z,W)\leq \psimod(Y,W).$$
\end{itemize}
\end{lem}
*****************************************************}

\subsubsection{Boundary of puzzle pieces}\label{boundary}
The ideal boundary of a $\psi$-puzzle piece $\BY$
is tiled by (finitely many) arcs $\lambda_i\subset \di \hat A_i$
that cover the ideal boundary of $\BU^m$ (where $m=\depth Y$) 
and arcs $\xi_i, \eta_i \subset \di \hat A_i$ mapped onto the Julia set $J(f)$. 
The arc $\la_i$ meets each $\xi_i,\, \eta_i$ at a single boundary point
corresponding to a path $\de: [0,1)\mapsto A$ that wraps around $K(f)$ infinitely many times,
while $\eta_i$ meets $\xi_{i+1}$ at a vertex $v_i\in Y\cap K(f)$. 
We say  that the arcs $\la_i$ form the {\it outer boundary } (or ``$O$-boundary'')
 $\di_O\BY$ of the puzzle piece $\BY$,
while the arcs $\xi_i$ and $\eta_i$
form its {\it $J$-boundary} $\di_J \BY$. 
Given a vertex $v=v_i$ of a puzzle piece $Y$, let $\di^v\BY=\eta_i\cup \xi_{i+1}$ 
stand for the part of the $J$-boundary of $\BY$ attached to $v$.

\section{Transfer Principle}
\label{covering lemma sec}

Let us now formulate two analytic results
which will play a crucial role in what follows.
The first one appeares in \S 2.10.3 of \cite{covering lemma}:

\newtheorem*{QA}{Quasi-Additivity Law}
\begin{QA}
    Fix some $\eta\in (0,1)$.
Let $\BV$ be a topological disk,  let  $K_i\Subset \BV $, $i=1,\dots, m$,
be pairwise disjoint full compact continua, 
and let $\phi_i: \A(1,r_i)\ra \BV \sm   \cup K_j$ 
be holomorphic annuli such that each $\phi_i$ is
an  embedding of some proper collar of $\T$  to a proper collar of $\di K_i$.
Then there exists a $\de_0 > 0$ (depending on $\eta$ and $m$) such that: \\
If  for some $\de\in (0, \de_0)$,
$\mod (\BV, K_i) < \de$ while $\log r_i > 2\pi \eta\de$  for all $i$,
then
$$
  \mod(\BV,  \cup K_i ) < \frac{2 \eta^{-1} \de}{m}.
$$
\end{QA}

The next result appears in \S 3.1.5 of \cite{covering lemma}:

\newtheorem*{Covering}{Covering Lemma}
\begin{Covering}
Fix some  $\eta\in (0,1)$. Let us consider
two topological disks $\BU$ and $\BV$, two full continua $A'\subset \BU$ and $B'\subset \BV$,
and two full compact continua $A\Subset A'$ and $B\Subset B'$.

  Let $f:\BU \ra \BV$ be a branched covering of degree $D$ such that
$A'$ is a component of $f^{-1}(B')$, and $A$ is the union of some components of $f^{-1}(B)$.
Let $d=\deg(f : A'\ra B')$.\\
Let $B'$ be also embedded into another topological disk $\BB'$.
Assume $\BB'$ is immersed into
$\BV$ by a map $i$  in such a way that $i|\, B'=\id$, $i^{-1}(B')=B'$,
 and $i(\BB') \sm B'$ does not contain the critical values of $f$.  \\
Under the following  ``Collar Assumption'':
$$
   \mod(\BB', B) > \eta \mod(\BU,  A),
$$
if 
$$ \mod(\BU, A) < \eps(\eta, D) $$
 then
$$
      \mod (\BV, B)  <  2 \eta^{-1} d^2 \mod(\BU,  A). 
$$
\end{Covering}  

  We will now apply these two geometric results to a dynamical situation.
  Recall from \S \ref{modified nest} that $\chi$ stands for the height of the Principal Nest, 
so that the quadratic-like map $g_\chi: E^\chi\ra E^{\chi-1}$ represents the
renormalization of $f$ with the filled Julia set $\KK$.

\begin{figure}[htbp]
\begin{center}
{\makeabbrevs
\input{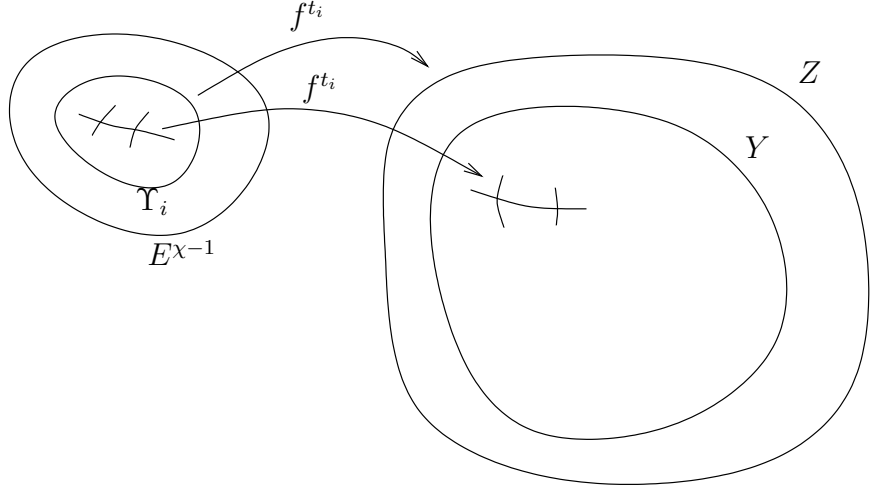}
\caption{The Transfer Principle}
\label{transfer}}
\end{center}
\end{figure}

\newtheorem*{Transfer}{Transfer Principle}
\begin{Transfer}
  Suppose there are two geometric puzzle pieces 
$Y\Subset Z$ with $\depth Z< \depth E^{\chi-1}$,  
and a sequence of moments of time $0<t_1<t_2<\dots <t_m$ such that:
\begin{itemize}
\item $t_m-t_1< p$ and $t_m<2p$;
\item   $f^{t_i}(\KK)\subset Y$;
\item $\Ups_i = f^{-t_i}(Z)|\, \KK \subset E^{\chi-1}$;
\item $\deg (f^{t_i}: \Ups_i\ra Z)\leq D$.
\end{itemize}
  Then there exist an absolute constant $C$ and $\eps=\eps(D)$ such that
$$
   \psimod(Z,Y)<  \frac{C}{m}\, \psimod(E^\chi, \KK),
$$
provided $\psimod(E^\chi, \KK) < \eps$.  
\end{Transfer}

\begin{pf}
Let $\KK_j=f^{t_j}(\KK)$.
We want to apply the Covering Lemma to the maps $f^{t_j}: (\Ups_j, \KK) \ra (Z,\KK_j )$. 
As the buffer around $\KK$ we take $E^{\chi+1}$.
Since
$$
   \depth E^{\chi+1}= \depth E^{\chi-1}+2p >\depth Z + t_j = \depth \Ups_j,
$$
we have $E^{\chi+1}\subset \Ups_j$ for any $j=1,\dots, m$. 

We let $\Om_j= f^{t_j}(E^{\chi+1})$
be the corresponding buffer around $\KK_j$.
Then $\deg(f^{t_j}: E^{\chi+1}\ra \Om_j)\leq 4$ since $t_j<2p$.
Moreover, 
\begin{equation}\label{bla-bla}
   \psimod(\Om_j, K_j) \geq  2\, \psimod (E^{\chi+1}, \KK) = \frac 1{2} \psimod (E^{\chi-1}, \KK) 
          \geq \frac 1{2} \psimod (\Ups_j, \KK),    
\end{equation}
which puts us in the position to apply the Covering Lemma with $\eta=1/2$
and $d=4$. It yields: 
\begin{equation}\label{comparable mod}
  \psimod(Z, \KK_j) \leq 2^6\, \psimod (E^{\chi-1}, \KK)= 2^7 \, \psimod (E^\chi, \KK),
\end{equation}
provided $\psimod(E^\chi,\KK)<\eps(D)$. 

Let us define $k_j$ as  $t_j$ if $t_j<p$,
and as $t_j-p$ otherwise. 
Then $k_j$'s are pairwise different numbers in between $0$ and $p$,
and hence the sets  $f^{k_j}(E^\chi)$ are pairwise disjoint.
Since $\Om_j\subset f^{k_j}(E^\chi)$, the buffers $\Om_j$ are pairwise disjoint as well.  
Moreover, by (\ref{bla-bla})
$$
  \psimod(\Om_j, \KK_j)\geq  \psimod(E^\chi, \KK),
$$
which, together with (\ref{comparable mod}), 
puts us into a position to apply the Quasi-Additivity Law with $\eta=2^{-7}$. 
It yields
$$
   \psimod(Z, Y) \leq \frac{2^{15}}{m} \psimod(E^\chi, \KK), 
$$
provided $\psimod(E^\chi,\KK)<\eps(D)$.
\end{pf}

\section{Improving the moduli}\label{moduli}

In this section we will prove Theorem \ref{high periods thm} for $\psi$-ql maps. 

Let $f_i: \BU_i' \ra \BU_i$ be a sequence of renormalizable $\psi$-ql maps
satisfying the $\eta$-molecule condition.
Let $p_i\to \infty$ stand for the renormalization periods of the $f_i$,
and let $\psimod (R f_i)\to 0$.
We need to show that 
$$
        \psimod (R f_i)/ \psimod (f_i)\to \infty.
$$

Let $P_{c_i}: z\mapsto z^2 + c_i$ be the straigtenings of the $f_i$.
Without loss of generality we can assume that $c_i\to c$.
Then the $\eta$-molecule condition implies 
that the quadratic polynomial $P_c$ satisfies the $(\bar\r, \bar\q, \bar\n)$-condition,
with  $\bar\r$, $\bar\q$ and $\bar\n$ depending only on $\eta $.    
Hence all nearby maps satisfy the $(\bar\r, \bar\q, \bar\n)$-condition as well.  
%
In what follows, we will fix one of these maps, $f=f_i$, with parameters $(\r,\q,\n)\leq (\bar\r, \bar\q, \bar\n)$, 
and consider its puzzle as described in \S \ref{principal nest}.   
All the objects under consideration (e.g., the principal nest $E^0\supset E^1\supset\dots$)
will be assoiciated with $f$ without making it notationally explicit.


\subsection{From the bottom to the top of the Principal Nest}

The following result proved  \cite{decorations} (Lemma~5.3)  
shows that if $\chi$ is big while the modulus $\psimod(E^{\chi-1}, E^\chi)$ 
is small then  $\psimod(E^0, E^1)$ is even smaller:

\begin{lem}\label{exp decay}
   For any $\k\in\N$ and  $\rho\in (0,1)$, there exist $\eps>0$ and $\underline\chi\in \N$
such that if $\chi\geq \underline\chi$ and $\psimod(E^{\chi-1}, E^\chi)<\eps$,
then 
$$
   \psimod(E^0, E^1)<  \rho\, \psimod(E^{\chi-1}, E^\chi) \leq  \rho\, \mod(\BE^\chi, \KK).
$$  
\end{lem}

 \begin{cor}\label{l is big - 3}
   For any $\k\in \N$ and $\rho\in (0,1)$, there exist $\eps>0$ and $\underline \chi\in \N$
such that if  $\chi\geq \underline{\chi}$, 
and $\mod(\BE^{\chi-1}, \KK)<\eps$, then for some puzzle piece $Y^m(z)$ we have:  
$$
   \psimod(Y^0(z), Y^m (z))<  \rho\, \mod(\BE^\chi, \KK).
$$  
\end{cor}

\begin{pf}
Let us apply $f^{\r(\q\n+1)}$ to the pair $(E^0, E^1)$. 
It maps $E^1$ onto some puzzle piece $Y^m (z)$, and maps $E^0$ onto $Y^0(z)$ with degree at most 
$2^{\r(\q\n+1)}$. 
By Lemma \ref{psi-maps}, 
$$
     \psimod( Y^0(z), Y^m (z))  \leq  2^{\r(\q\n+1)}  \psimod ( E^0, E^1).
$$
Together with Lemma \ref{exp decay} this yields the assertion.
\end{pf}

\subsection{Around the stars}
  
We will now go back to the original map $f$. 
Recall that $p$ stands for its renormalization period,
$\k=\r\q\n$, and $\bl$ is introduced in Lemma \ref{single point}. 

\begin{lem}\label{p is big}
   For any $\k\in \N$ and $\rho\in (0,1)$, there exist $\eps>0$ and $\underline p\in \N$
such that if  $p\geq \underline{p}$, 
and $\mod(\BE^{\chi-1}, K)<\eps$, 
then either for some puzzle piece $Y^{\bl}(z)$,  
$$
    0 <  \psimod(Y^0(z), Y^{\bl} (z))<  \rho\, \mod(\BE^\chi, \KK),
$$  
or for some periodic point $\alpha_\mu\in \balpha$, 
$$
    0< \psimod(S^1(\alpha_\mu), S^{\bl}(\alpha_\mu)) <  \rho\, \mod(\BE^\chi, \KK).
$$

\end{lem}

\begin{pf}
  By Corollary \ref{l is big - 3}, it is true when  $\chi\geq \underline{\chi}$,
so assume  $\chi\leq \underline{\chi}$. 
It will follow from the Trnasfer Principle of \S \ref{covering lemma sec}. 
Let $E^{\chi-1} =  Y^{\tau_0}$.
Note that
$$
    \deg (f^{\tau_0}|\, E^{\chi-1} ) \leq 2^{\underline{\chi}+\r(\q\n+1)},
$$    
so it is bounded it terms of $\k$ and $\rho$. 

Let us then select the first moment $\tau\geq \tau_0$ such that $f^\tau(K)\not\subset S^{\bl}(\balpha)$.
It is bounded by $p+\k$ (since  $f^{p+\k}(K) = f^{\k} (K)\not\subset S^{\bl}(\balpha) $).   

Let   $m=[\frac{C}{\rho}]+1$, where $C$ is the constant from the Transfer Prinsiple.
Let $s=s(\r,\q)$ be the number of puzzle pieces $Y^{\bl}(z)$ 
in the complement of the star $S^{\bl}(\balpha)$ (see \S \ref{stars}),
 and let  $N= \k  s m^2 $.
Let us consider the piece of $\orb \KK$ of length $N$,    
\begin{equation}\label{orbit}
      f^t(\KK), \quad t=\tau, \dots, \tau+N. 
\end{equation}
Then one of the following options takes place:

\ssk\nin (i) $m$ sets $\KK_j= f^{t_j}\KK$, $j=0,\dots, m-1$,  in the orbit (\ref{orbit})
  belong to some puzzle piece $Y^{\bl}(z)$ 
in the complement of the star $S^{\bl}(\balpha)$;

\ssk\nin (ii)  $\k  m $ consecutive sets $f^t(\KK)$, $t=i+1,\dots i+\k m$,
 in the orbit (\ref{orbit})
belong to the star $S^{\bl}(\balpha)$.
Here $i$ is selected so that that $f^i(\KK)$ does not belong to the star $S^{\bl}(\balpha)$.
Note that $i\geq\tau\geq\tau_0$ by definition of $\tau$. 

\msk
Assume the first option occured. 
Then let us consider the puzzle piece $Y^0(z)$  of depth $0$ containing $Y^{\bl}(z)$.
By Lemma \ref{compactly conained}, $Y^{\bl}(z)\Subset Y^0(z)$.
Let us consider if pullback $\Ups_j= f^{-t_j}(Y^0(z))|\KK$ containing $\KK$. 
Then $\Ups_j\subset E^{\chi-1}$ (since $t_j\geq \tau_0)$). 

Moreover,  the map $g^{t_j} : \Ups_j\ra Y^0(z)$ has degree bounded in terms of $\k$ and $\rho$.  
Indeed,   $\deg (f^{\tau_0}|\, E^\chi )$ and $N$, $\bl$, $\k$ are bounded in these terms.
Hence it is enough to show that the trajectory   
$$
   f^m (\Ups_j), \quad \tau_0<  m < \tau-\bl-\k,
$$
does not hit the critical point.
But if $f^m(\Ups_j)\ni 0$ then  $f^{m+\k}(\Ups_j)$ would  land 
outside $S^{\bl}(\balpha)$ (since $f^\k(0)\not \in S^\bl(\balpha)$ and $\depth f^{m+\k}(\Ups_j)\geq \bl$).
Then $f^{m+\k}(\KK)$ would land outside $S^{\bl}(\balpha)$ as well
contradicting the defintion of $\tau$ as the first landing moment 
of $\orb K$ in $S^\bl(\balpha)$ after $\tau_0$. 

Now,  selecting $\underline{p}$ bigger than $\k+N$,
we bring ourselves in the position to apply the Transfer Principle
with $Y=Y^\bl(z)$, $Z=Y^0(z)$. 
It yields:
$$
   \psimod(Y^0(z), Y^{\bl}(z) ) \leq \frac{C}{m} \psimod(E^\chi, \KK)\leq \rho\,  \psimod(E^\chi, \KK).
$$

\msk Assume now the second option occured.
Then there is a point $\alpha_\mu\in \balpha$ such that  
$f^{i+\k j}(\KK)\subset S^\bl(\alpha_\mu)$ for  $j=1,\dots, m$,
while $f^{i}(\KK)\subset S^\bl(\alpha_\mu')$.
Let us pull the star $S^1(\alpha_\mu)$ back by $f^{i+\k j}$:   
$$
  \Ups_j= f^{-(i+\k j)}(S^1(\alpha_\mu))|\, \KK,\quad  j=0,\dots, m-1.
$$ 
Let us show that 
\begin{equation}\label{inclusion}
   \Ups_j\subset E^{\chi-1}. 
\end{equation}
We fix some $j$ and let $\Ups=\Ups_j$. 
%
We claim  that $\inter f^i(\Ups)$ does not contain any points of $\balpha$. 
Since $\alpha_{\mu}$ is the only point of $\balpha$  inside $\inter S^1(\alpha_\mu)$, 
it is the only point of $\balpha$ that can be inside $\inter f^i(\Ups)$. 
If $\mu=0$ then $f^i (\KK)\subset S^\bl(\alpha')$, 
so by Lemma~\ref{qn-pullbacks}, 
$f^i(\Ups) \subset S^1(\alpha')$, 
which does not contain $\alpha$. 
If $\mu\not=0$ then the points $\alpha_\mu$ and $\alpha_\mu'$ are separated by $\alpha$
in the filled Julia set. Since $f^i(\Ups)$ is a geometric puzzle piece containing both 
$\alpha_\mu$ and $\alpha_\mu'$, it must contain $\alpha$ as well --  contradiction.  
The claim follows.

Thus,  $f^i(\Ups)$ is a geometric puzzle piece whose interior does not contain any points of $\balpha$.
Hence it is contained in some puzzle piece $Y^0(z)$ of zero depth. 
Then $\Ups \subset f^{-i}(Y^0(z))|\, 0=Y^i$.
But since $i \geq \tau_0$,  
$
   Y^i  \subset Y^{\tau_0} =  E^{\chi-1},
$
and (\ref{inclusion}) follows. 

Other assumptions of the Transfer Principle 
(with $Y=S^\bl(\alpha_\mu)$  and $Z=S^1(\alpha_\mu)$) 
are valid for the same reason as in the first case.
The lemma follows.
\end{pf}

\subsection{Bigons}

 A geometric puzzle piece with two vertices is called a {\it bigon},
and the corresponding pseudo-puzzle piece is called a $\psi$-bigon. 
Given  a bigon $Y$ with vertices $v$ and $w$, 
let $S_Y^v$ and $S_Y^w$ stand for the components of $K(f)\sm \inter Y$ containing $v$ and $w$ 
respectively.  We let $S_Y=S^v_Y\cup S^w_Y$.

Recall from \S \ref{boundary} 
that the ideal boundary of the corresponding $\psi$-puzzle piece $\BY$ comprises the outer boundery  
$\di_O \BY$ (in the case of bigon consisting of two arcs) and the $J$-boundary
$\di_J\BY= \di^v\BY\cup \di^w\BY$ attached to the vertices. 
Let $\GG_\BY=\GG_\BY(v,w)$ stand (in the case of bigon)                    
for the family of horizontal  curves in $\BY$ connecting $\di^v\BY$ to $\di^w\BY$,
and let $\Bd_\BY(v,w)$ stand for its extremal length.

More generally, let us consider a puzzle piece $Y$ whose vertices are bi-colored, 
i.e., they are partitioned into two non-empty subsets, $B$ and $W$. 
This induces a natural bi-coloring of $K(f)\sm \inter Y$ and of $\di_J\BY$:
namely, a component of these sets  attached to a black/white vertex inherits the
corresponding color. 
Let $\GG_\BY$ stands for the family of horizontal  curves in $\BY$ 
connecting boundary components with different colors. 

For a geometric puzzle piece $Y$, 
let $v(Y) \subset K$ denote the set of vertices of $Y$.
Suppose that $Y \subset Z$ are (geometric) puzzle pieces with the same equipotential depth;
we say that $Y$ is \emph{cut out of} $Z$ if $v(Z) \subset v(Y)$,
so that we have produced $Y$ by cutting out pieces of $Z$.
If the vertices of $Y$ are bicolored,
then the vertices of $Z$ are as well.

\begin{lem} \label{cut out}
  If $Y$ is cut out of $Z$,
  and $v(Y)$ is bicolored,
  then the family of curves $\GG_Z$ overflows $\GG_Y$.
\end{lem}

\begin{pf}
  We prove the Lemma by induction on the cardinality of $v(Y) \sm v(Z)$.
  First suppose that $v(Y) = v(Z) \cup \{w\}$, 
  where $w \notin v(Z)$.
  Let $\gamma \in \GG_Z$;
  the two endpoints of $\gamma$ lie in differently-colored components $\di^x \BZ$, $\di^y \BZ$
  of $\di_J \BZ$.
  If $\gamma$ lifts to $\GG_Y$, 
  then we are finished.
  Otherwise, 
  we can start lifting $\gamma$ from the endpoint (of $\gamma$)
  that lies in the component (say $\di^x Z$) of $\di_J Z$
  whose color is different from that of $w$.
  Then that partial lift of $\gamma$ will connect $\di^x Y$ and $\di^w Y$
  and hence will belong to $\GG_Y$.
  
  In general,
  if $|v(Y)| > |v(Z)| + 1$,
  we can let $Y'$ be such that 
  $v(Y) = v(Y') \cup \{w\}$,
  and $v(Y') \supset v(W)$.
  Then given $\gamma \in \GG_W$,
  we can lift part of $\gamma$ to $\GG_{Y'}$ by induction,
  and then to $\GG_Y$ as before.
\end{pf}

Given a puzzle piece $Y$ with $v(Y)$ bicolored
and  a bigon $P$,
we let $P\succ Y$ if the vertices of $P$ belong to components of $K(f)\sm \inter Y$
with different colors,
while the equipotential depths of $P$ and $Y$ are the same. 

\begin{lem}\label{overflowing}
  If $P\succ Y$ then the family $\GG_P$ overflows  $\GG_Y$.
\end{lem}

\begin{pf}
  Let $P'$ be the bigon whose vertices are the vertices of $Y$ that separate $\inter Y$
  from the $v(P)$.
  Then $\GG_{P}$ overflows $\GG_{P'}$,
  and $\GG_{P'}$ overflows $\GG_Y$ by Lemma \ref{overflowing}.
\end{pf}

Let $\WW_Y$ stand for the width of $\GG_\BY$. 

\begin{lem}\label{WW-est}
  Let $Y$ and $P$ be two bigons such that the vertices of $f^n(P)$  are separated by $Y$ for some $n$,
and the equipotential depth of $Y$ is $2^n$ times bigger than the equipotential depth of $P$.
Then $\WW_Y \geq 2^{-n} \WW_P$.  
\end{lem}

\begin{pf}
  Since the vertices of $f^n(P)$ are separted by $Y$, 
there is a component $Z$ of $f^{-n}(Y)$
such that $P\succ Z$. 
By Lemma \ref{overflowing}, $\WW_P\leq \WW_Z$.   
On the other hand, the map $f^n: Z\ra Y$ has degree at most $2^n$,
and maps horizontal curves in $Z$ to horizontal curves in $Y$.     
Hence $\WW_Z\leq 2^n \WW_Y$. 
\end{pf}

\begin{lem}\label{comb separation of vertices}
   Let $Y$ be a bigon with vertices $u$ and $v$ of depths $l$ and $m$ 
satisfying the following property:
If $l=m$ then $f^l u\not=f^l v$. Then there exists an $n\leq \max(l,m)+\r$ such that
$f^n u$ and $f^n v$ are separated by the puzzle piece $Y^1$. 
\end{lem}
 
\begin{pf}
  By symmetry, we can assume that $l \le m$. 
  Suppose that $f^m u \neq f^m v$;
  then we can find $0 \le t < r$ such that 
  $f^{m+t} u$ and $f^{m+t} v$ are on opposite sides of the critical point,
  so they are separated by $Y^1$.
  Otherwise, 
  we must have $l < m$,
  and then $f^{m-1} u = - f^{m-1} v$;
  hence they are separated by $Y^1$.
\end{pf}

\subsection{Amplification}
We can now put together all the above results of this section as follows:

\begin{lem}\label{amplification}
 For any $\k\in \N$ and $\rho\in (0,1)$, there exist $\eps>0$ and $\underline p\in \N$
such that if  $p\geq \underline{p}$, 
and $\mod(\BE^{\chi-1}, \KK)<\eps$, then
$$
     \Bd_{Y^1}(\alpha, \alpha') \leq \rho \mod(\BE^{\chi-1}, \KK) . 
$$   
\end{lem}

\begin{pf}
  Under our circumstances,
  Corollary \ref{l is big - 3} and Lemma \ref{p is big} imply that there exist geometric puzzle pieces $Y\Subset Z$ 
with the bidepth of $Z$  bounded by $(1,1)$ while the  
bidepth of $Y$  bounded in terms of $\k$,  such that
$$
   \psimod(Z,Y)\leq \rho  \mod(\BE^{\chi-1}, \KK).
$$
For any vertex $v$ of $Z$, there exists a vertex $v'$ of $Y$ such that the rays of $\di Y$ landing at $v'$
separate $\inter Y$ from $v$. These two rays together with the two rays landing at $v$ 
(trancated by the equipotential of $Z$) form a bigon $B^v$. By the Parallel Law,
there exists a vertex $v$ of $Z$ such that
$$
          \Bd_{B^v}(v,v')\leq N \psimod (Z,Y),
$$                                                                 
where $N\leq N(\r)$ is the number of vertices of $Z$.    
By Lemma \ref{comb separation of vertices},
there is an iterate $f^n(B^v)$ such that the vertices $f^n(v)$ and $f^n(v')$ are separated by 
$\inter Y^1$. By Lemma \ref{WW-est}, 
$$
        \Bd_{Y^1} (\alpha, \alpha') \leq 2^n  \Bd_{B^v}(v,v').
$$
Putting  the above three estimates together, we obtain the assertion.     
\end{pf}

\begin{rem} The name ``amplification'' alludes to the extremal {\it width} which is amplified 
  under the push-forward procedure described above.
\end{rem} 

\subsection{Separation}

  The final step of our argument is to show that the vertices $\alpha$ and $\alpha'$ are
well separated in the bigon $Y^1$. 

\comm{************
   In this section we will show that   the modulus $\Bd_{Y^1}(\alpha, \alpha')$
that measures  the extremal distance between $L$ and $R$ is comparable
with $\mu=\mod(\BU, K)$. 

Let $Y$ be a geometric puzzle piece.  For  two vertices $v$ and $w$ of $Y$,
 we let  
$$
 \BW_Y(v,w) =\WW(\GG_Y^{vw}).
$$
We define the {\it pseudo-conductance} of $Y$ as 
$$
    \BW_Y =\sup_{v,w} \BW_Y(v,w),
$$
where the supremum is taken over all pairs of the vertices of $Y$.        
*******************}

\begin{lem}\label{main separation} 
 There exists $\kappa=\kappa(\r,\q,\n)>0$ such that   
   $$   \Bd_{Y^1}(\alpha, \alpha') \geq \kappa \mod (\BU,K). $$
\end{lem}

\nin {\it Idea of the proof.}
The proof is the same as the one of Proposition 5.12 in \cite{decorations},
so we will only give an idea here. 

Let $\BY$ be a $\psi$-puzzle piece, and let $v$ and $w$ be two vertices of it.
A {\it multicurve in $\BY$ connecting $\di^v\BY$ to $\di^w\BY$} is a sequence of proper paths
$\gamma_i$, $i=1,\dots, n$,  in $\BY$ connecting $\di^{v_{i-1}}\BY$ to $\di^{v_i}\BY$,
where $v=v_0, v_1,\dots, v_n=w$ is a sequence of vertices in $\BY$.  
Let $\BW_\BY(v, w)$ stand for the extremal width of the family of multicurves in $\BY$ 
connecting $\di^v\BY$ to $\di^w\BY$.
Let 
$$
    \BW_Y =\sup_{v,w} \BW_Y(v,w). 
$$

Let us estimate this width for the puzzle piece $P$ introduced in \S \ref{P-sec}.
To this end let us consider puzzle pieces $Q^v$ from Lemma \ref{Qv}. 
Let $r$ be the depth of these puzzle pieces,
$T^{vw}=\cl(K_P\sm (Q^v\cup Q^w))$, and let $v'=Q^v\cap T^{vw}$, $w'=Q^w\cap T^{vw}$. 
For any multicurve $\gamma$ in  $\BP$ connecting $\di^v \BP$ to $\di^w \BP$, 
one of the following events can happen:

\ssk\nin (i) $\gamma$ skips over $T^{vw}$;

\ssk\nin (ii) $\gamma$ contains an arc $\gamma'$ connecting an equipotential of depth $r$ 
             to $T^{vw}$;

\ssk\nin (iii)  $\gamma$ contains two disjoint multicurves,
$\de^v$ and $\de^w$, that do not cross this equipotential  and such that 
$\de^v$ connects $\di^v\BQ^v$ to $\di^{v'}\BQ^v$, while $\de^w$ connects $\di^{w'}\BQ^w$ to $\di^w\BQ^w$. 

\ssk
It is not hard to show that the width of the first two families of multicurves is $O(\mod(\BU,K))$
(see \S\S 5.5-5.6 of \cite{decorations} ). Concerning each family of multicurves $\de^v$ or $\de^w$
that appear  in (iii), it is conformally equivalent to a family of multicurves connecting appropriate
two vertices of $P$ (since $Q^v$ and $Q^w$ are conformal copies of $P$). By the Series and Parallel Laws,
$$
     \BW_P \leq \frac 1{2} \BW_P +O(\mod(\BU, K)),
$$
which implies the desired estimate.

\subsection{Conclusion}\label{conclusion}

Everything is now prepared for the main results.
 Lemmas \ref{amplification} and \ref{main separation} imply:

\begin{thm}[Improving of the moduli]\label{mod improve}
For any parameters $\bar\r,\bar\q,\bar\n$ and any $\rho>0$,
 there exist  ${\underline p} \in \N$ and $\eps>0$ with the following property.
Let $f$ be a renormalizable $\psi$-quadratic-like map with renormalization period $p$
satisfying the $(\bar\r, \bar\q, \bar\n)$-molecule condition, 
and let $g$ be its first renormalization.
Then 
$$
 \{ p\geq \underline p\ \ \mathrm{and}\ \psimod(g)< \eps \}  \imply \psimod (f) < \rho \mod(g).
$$
\end{thm}

  Theorem \ref{mod improve},
together with Lemma \ref{ql extension}, 
implies Theorem \ref{high periods thm} from the Introduction. 
The Main Theorem  follows from Theorem \ref{mod improve} combined
with the following result  (Theorem 9.1 from \cite{K}): 

\begin{thm}[Improving of the moduli: bounded period]\label{improving mod-2}
For any $\rho\in (0,1)$, there exists $\underline p=\underline p (\rho)$ 
such that for any $\bar p\geq \underline p$, there exists $\eps=\eps(\bar p)>0 $ with the following property.
Let $f$ be primitively renormalizable $\psi$-quadratic-like map, 
and let $g$ be the corresponding renormalization.
Then
$$
   \{\underline p\leq p\leq \bar p\ \ \mathrm{and}\ \psimod(g) < \eps\}   \imply \psimod(f) < \rho\, \psimod(g).
$$
\end{thm}

Putting the above two theorems together, we obtain:

\begin{cor}\label{improving-3}
 For any $(\bar\r, \bar\q, \bar\n)$,   there exist an $\eps>0$ and $l\in \N$ with the following property.
For any  infinitely renormalizable $\psi$-ql map $f$ satisfying the  $(\bar\r, \bar\q, \bar \n)$-molecule
condition
with renormalizations $g_n=R^n f$,
if $\psimod(g_n)< \eps$, $n\geq l$,  then $\psimod (g_{n-l})<\frac 1{2} \psimod(g_n)$.
\end{cor}

This implies the Main Theorem, in an important refined version.
We say that a family $\MM$ of little Mandelbrot copies (and the corresponding
renormalization combinatorics) has {\it beau}%
\footnote{According to Dennis Sullivan, 
   ``beau'' stands for ``bounded and eventually universal''. }
{\it a priori} bounds
if there exists an $\eps =\eps(\MM)>0$ and a function $N: \R_+\ra \N$
with the following property.
Let $f: U\ra V$ be a quadratic-like map with $\mod(V\sm U)\geq \de>0$
that is  at least $N=N(\de)$ times renormalizable. Then for any $n\geq N$,
the $n$-fold renormalization of $f$ can be represented  by a quadratic-like map $R^n f: U_n\ra V_n$ with
$\mod(V_n\sm U_n)\geq \eps$.

\newtheorem*{Beau}{Beau Bounds}  
\begin{Beau}
  For any parameters $(\bar\r, \bar\q, \bar\n)$, 
the family of renormalization combinatorics satisfying the  $(\bar\r, \bar\q, \bar\n)$-molecule condition has
beau a priori bounds.
\end{Beau}

\subsection{Table of notations}

$p$ is the renormalization period of $f$;

$\gamma$ is a dividing periodic point of period $\bt$, 

$\bgamma$ is its cycle;

$\s$ is the numner of rays landing at $\gamma$;   

$\alpha$ is a dividing periodic point of periods $\r=\bt\s$;

$\alpha'=-\alpha$, $\alpha_j=f^j\alpha$;

$\balpha=\{\alpha_j\}_{j=0}^{p-1}$ is the cycle of $\alpha$; 

$\q$ is the number of rays landing at $\alpha$;

$\n$ is the first moment such that $f^{\r\q\n}(0)$ is separated from $0$ by the rays 
      landing at $\alpha'$,

 $\k= \r\q\n$;

$\bl= \k+1+\r^2$ is such a depth that the stars $S^\bl(\alpha_j)$, $j=0,1,\dots, \r-1$,
       are all disjoint;

\end{document}